# Soñando con números, María Andresa Casamayor (1720-1780)

Julio Bernués y Pedro J. Miana

**Resumen:** La zaragozana María Andresa Casamayor es conocida por ser la primera mujer escritora de un texto científico en España. En este artículo respondemos a varios de los enigmas más importantes de su desconocida biografía, como la fecha de su nacimiento, el origen de su familia, la localización de la casa en la que habitó así como las de sus familiares, su profesión como maestra de primeras letras e incluso… su verdadero nombre.

---

Hace casi 300 años nacía en la Zaragoza del Siglo de las Luces María Andresa Casamayor de La Coma. Con poco más de 17 años, publicó su obra *Tyrocinio arithmetico, Instrucción de las quatro reglas llanas.* Esta modesta e interesante obra fue el primer libro de ciencias publicado por una mujer en España. Su autora, bastante olvidada, encarna por su vocación y dedicación el ejemplo de mujer científica que desarrolla su trabajo en circunstancias tremendamente desfavorables.

**Comercio y educación en la Zaragoza del siglo XVIII**

La colonia francesa de Zaragoza dominaba el comercio en Aragón entre los reinos de España y el resto de Europa formando a principios del siglo XVIII un numeroso grupo de población con fuertes interrelaciones comerciales y familiares (Salas, 2003). En este ambiente se sitúan los progenitores de nuestra protagonista, Juan Joseph Casamayor, un comerciante textil francés que casará en 1705[1] con Juana Rosa de La Coma, hija de comerciantes de ascendencia también francesa pero ya afincados en Zaragoza. El matrimonio tendrá 9 hijos e hijas[2] y entre éllos, *María* Juana Rosa *Andresa* nacerá un 30 de noviembre, día de San Andrés, de 1720 siendo bautizada al día siguiente en la iglesia del Pilar.

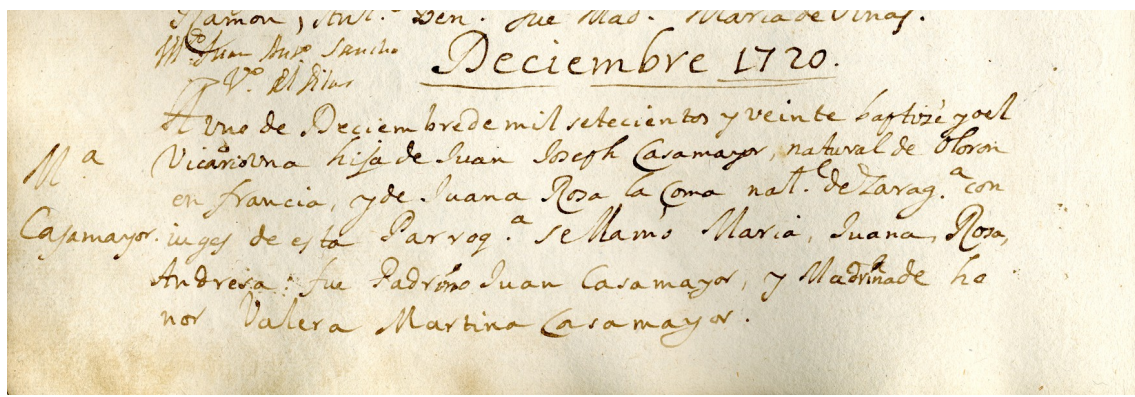

Apunte del bautismo de María Andresa Casamayor de La Coma

Es muy probable que, como hacían muchas otras familias de Zaragoza acomodadas como la suya, María Andresa recibiera sus *primeras letras* de forma colectiva con sus hermanos y hermanas en la casa en la que la familia vivía alquilada sita en la Calle del Pilar[3]. En este ambiente, la pequeña María Andresa muy pronto destacó.

Leer, escribir, contar, operar con las cuatro reglas… eran habilidades de obligado conocimiento para el comercio, la principal actividad que se realizaba en el entorno de María Andresa, entre los diversos territorios aragoneses, castellanos y franceses con un amplio número de unidades de moneda, longitud, superficie o peso.

Los ilustrados del siglo XVIII consideraron el mejorar la educación como uno de los objetivos fundamentales de sus políticas para el progreso del país. Y es que en los años de infancia de María Andresa, el porcentaje de analfabetos era enorme, especialmente entre las mujeres (Alfaro, 2017). Aunque la implicación pública fue en aumento sobre todo a partir del último tercio del siglo, la enseñanza en el siglo XVIII era principalmente impartida en centros religiosos.

De especial importancia para nuestra historia es el papel de padres escolapios que llegan a Zaragoza el 27 de octubre de 1731. Los escolapios impartían una enseñanza de calidad, gratuita y universal (eso sí, sólo para hombres) tanto en letras como en ciencias. Uno de estos primeros escolapios fue Juan Francisco de Jesús, catedrático de Matemáticas y uno de los censores del libro del Tyrocinio. El colegio fundacional de Zaragoza, conocido hoy como "el colegio de Conde Aranda", se abrió el 19 de febrero de 1740 con el nombre de Colegio de Santo Tomás de Aquino de las Escuelas Pías de Zaragoza en homenaje al patrocinio del arzobispo Tomás Crespo Agüero. Este prelado, preocupado también por la educación femenina, encargó ésta a las "Señoras de la Real Casa de la Enseñanza".

Zaragoza contaba además con 10 maestros de primeras letras (para niños) que estaban distribuídos por la ciudad e impartían clases en locales o en sus propios domicilios. La organización de ese tipo de enseñanza para niñas será algo posterior (Domínguez, 1999). Por otro lado la ciudad disponía de una escuela pública de primeras letras, las *aulas públicas*, situada en el edificio jesuíta del Colegio del Padre Eterno.

Los políticos ilustrados realizaron numerosos esfuerzos para organizar la enseñanza pública y de hecho, a pesar de su precario funcionamiento, un alto porcentaje de los pueblos aragoneses contarán al acabar el siglo con la conocida figura del *"maestro de escuela"*, que será (muy mal) pagado bien por ayuntamientos o … por las propias familias de los alumnos (Alfaro, 2017. Domínguez, 1999).

**El Tyrocinio Arithmetico**

El libro del Tyrocinio ha sido hasta hoy la principal fuente de información sobre su joven autora. Publicado en marzo de 1738, es un libro de tamaño de cuarto, de 16 cm x 22 cm aproximadamente, de 78 páginas con algún error de paginación. Una copia digitalizada puede encontrarse en la web de la Biblioteca Nacional.

El cultismo latino *Tyrocinio* (que significa, aprendizaje o formación) era un vocablo utilizado en el siglo XVIII para titular obras de ciencias o de letras.

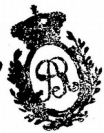

Portada del Tyrocinio Arithmetico

Desde el punto de vista matemático, el libro está escrito en un lenguaje ágil y eminentemente práctico, con una gran cantidad de ejemplos y casos reales que permiten al lector aprender de forma directa el manejo de las cuatro reglas del *álgebra menor*: suma, resta multiplicación y división. Además muestra un conocimiento preciso de las unidades que se manejaban a diario en el comercio de principios del siglo XVIII. Si a ello unimos que en la introducción se resalta el carácter pedagógico de la obra, en el perfil de María Andresa Casamayor podríamos destacar una gran habilidad aritmética y una profunda preocupación por la educación. En este sentido, la obra de María Andresa se adelanta en varias décadas a lo que será el modo de hacer de las mujeres ilustradas.

En el plano personal, la autora firma con un pseudónimo masculino, *Casandro Mamés de La Marca y Araioa*. Esta firma es un perfecto anagrama, misma letras en diferente orden, de su nombre **María Andresa Casamayor de La Coma**. Curiosamente esta noticia, sacada de la cita de Félix Latassa en su monumental obra *"Biblioteca nueva de los escritores aragoneses"*, contiene un error ya que Latassa llama a nuestra protagonista "María Andrea" en lugar de María Andresa, un equívoco que ha llegado hasta nuestros días.

*Casandro* se reconoce como "discípulo de la Escuela Pía" y dedica el libro a la misma "Escuela Pía del Colegio de Santo Thomás de Zaragoza". Esta dedicatoria parece significar un apoyo a los recién llegados escolapios.

Firma la dedicatoria en Almodóvar de el Pinar. Este pueblo de la serranía de Cuenca, cercano a Teruel, contaba con un Colegio Escolapio desde 1724. Fue conocido como el "pueblo de las carretas" al ser un importante nudo de comunicaciones en la época. Se desconoce el interés de la autora al nombrar esta localidad.

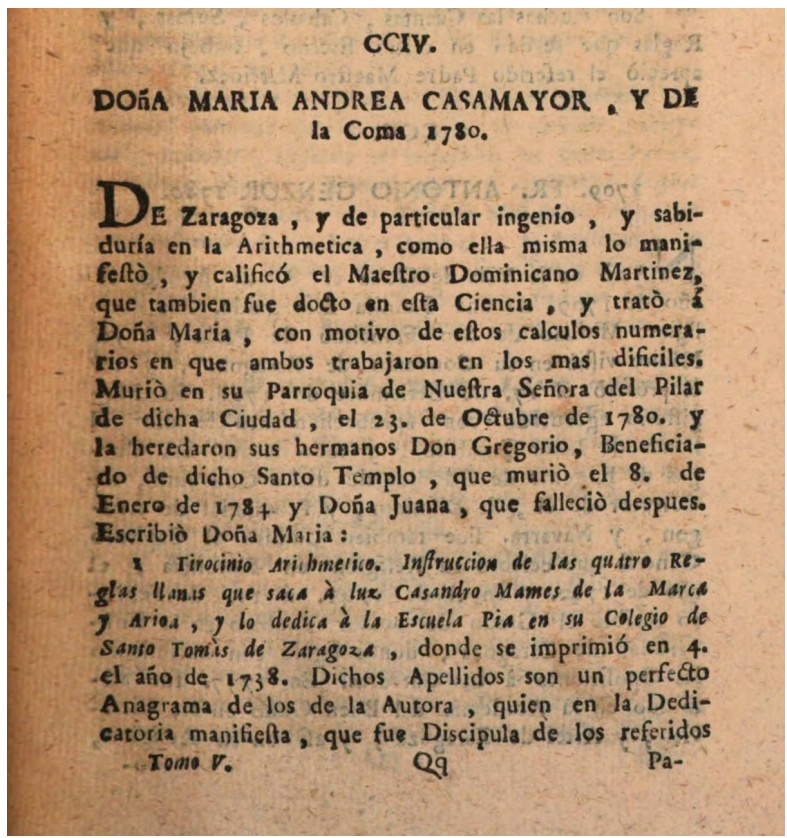

Apunte biográfico de "María Andrea" en la enciclopedia de Latassa

Siguiendo las normas de la época, el Juez de Impresiones y Oidor de la Real Audiencia, Don Alonso Pérez de Mena, solicita al escolapio y catedrático de matemáticas Juan Francisco de Jesús y al fraile Pedro Martínez, Regente de Estudios del Colegio de San Vicente Ferrer de la Orden de Predicadores, las censuras correspondientes (hoy también las llamaríamos reseñas) como trámite para aprobar la publicación del Tyrocinio. Pedro Martínez nos cuenta cómo otros libros similares a éste suelen ser más extensos, lo que encarece su coste, y nos desvela los motivos del autor (autora) para publicar tan breve obra:

> *su fin, en esta Obrilla solo es facilitar esta instruccion a muchos, que no pueden lograrla de otro modo.*

**Y después del Tyrocinio…**

La segunda obra de Maria Andresa Casamayor fue el manuscrito, hoy perdido, *"El Para si solo"* de Casandro Mamés de la Marca y Arioa. *Noticias especulativas, y prácticas de los Numeros, uso de las Tablas de las Raizes, y Reglas Generales para responder à algunas Demandas, que con dichas Tablas se resuelven sin la Algebra.* De 109 hojas de tamaño folio es de un nivel matemático superior al *Tyronicio.* Tal y como comenta Latassa (1802),

> *Son muchas las Cuentas, Calculos, Sumas, y Reglas que se dán en dicho Escrito, trabajo que apreciò el referido Padre Maestro Martinez.*

Uno de los principales apoyos con el que contó la joven María Andresa Casamayor fue el del dominico y censor de su libro, Pedro Martínez. Latassa (1802) escribe los apuntes biográficos de ambos señalando sus colaboraciones mutuas.

Desconocemos si, debido al alto nivel matemático alcanzado por María Andresa, ésta llegó a entrar en contacto con alguno de los matemáticos afincados en Zaragoza o con alguno de los notables libros de matemáticas que vieron la luz en las imprentas zaragozanas de la época: En 1723 se publica *Euclides Geometria Especulativa y Practica de los Planos, y Solidos* del militar Antonio José Deu y Abella. Un año más tarde se reedita en *Arithmetica practica muy útil y necessaria para todo genero de Tratantes y Mercaderes* del valenciano Geronimo Cortés. Francisco Xavier García, residente en Zaragoza, publica su *Arithmetica especulativa y practica y arte mayor o algebra mayor o algebra* en 1733. Merece la pena mencionar que anteriormente en Barcelona en 1698 el ingeniero militar aragonés Francisco Larrando de Mauleón publicaba sus *Elementos de Euclides*.

**María Andresa, maestra de niñas.**

Inmediatamente después de escribir el Tyrocinio, trágicos acontecimientos van a determinar el futuro de María Andresa: Juan Jose Casamayor, su padre, fallece el 14 de marzo de 1738 y Fray Pedro Martínez, su amigo y colaborador, el 14 de noviembre de 1739. Por otro lado su familia cercana, los La Coma y en concreto el heredero de la familia Joseph de La Coma, iniciará en 1740 un proceso de endeudamiento que terminará en 1748 con la pérdida de todas sus casas[4]. También el arzobispo Tomás Crespo fallecerá el 3 de marzo de 1742 y a mediados del siglo XVIII, el matemático Juan Francisco de Jesús se trasladará al Colegio de los Escolapios de Valencia. De repente, todos los apoyos que había tenido la joven María Andresa Casamayor han desaparecido en pocos años.

A diferencia de lo que era habitual para una mujer de la sociedad zaragozana, María Andresa ni se casará ni entrará en la Iglesia, así que el resto de su vida deberá trabajar para ganarse la vida [5]. Fue maestra de niñas y durante buena parte de su vida, maestra de primeras letras en las aulas públicas de la ciudad. La publicación del *Tyrocinio arithmético* pudo servirle de carta de presentación para su labor profesional. Como parte de su retribución, le será facilitada una casa donde vivir.

Por el censo de población de 1766 que se conserva en el Archivo Histórico Municipal de Zaragoza sabemos que María Andresa Casamayor vivía (sola) en una casa de la calle Palomar que hace esquina con la actualmente llamada calle de la Viola que va a la plaza de San Agustín, en la Parroquia de Santa María Magdalena. La casa, que todavía existe en la actualidad, era propiedad de Joseph Lasala, escribano real. La anotación del censo indica "*Andresa Casamayor*, no paga".

![Apunte en el censo de población de 1766]

Apunte en el censo de población de 1766.

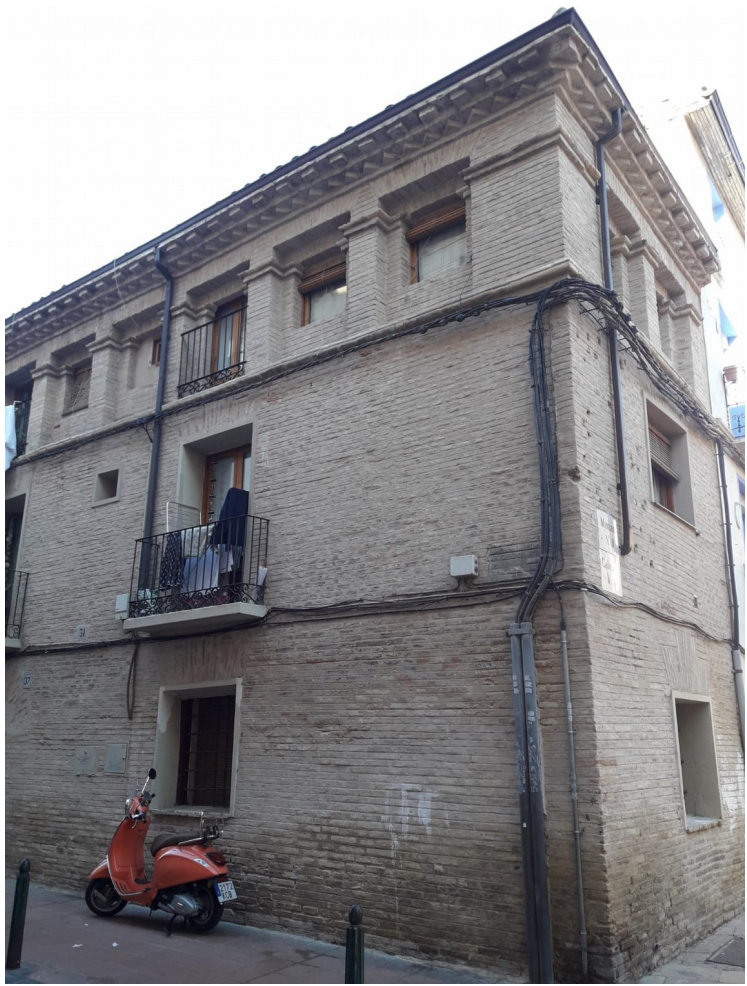

La casa de María Andresa en la actualidad

En el Archivo Municipal del Ayuntamiento de Zaragoza, se conserva el documento *Apuntación de las licencias que se han dado a las Maestras de Niñas para que puedan*

*enseñar,* que nos dice el lugar donde daba sus clases. En la primera línea del listado de maestras de niñas aparece *"Seminario Viejo. Maria Casamaior"* (Domínguez, 1999).

La historia de este lugar comienza en 1767 cuando se produce la expulsión de los Jesuítas y con élla la búsqueda de nuevos usos para sus instalaciones. Entre éllos el mencionado Colegio del Padre Eterno[6]. Este edificio había alojado, en su planta baja y desde al menos 1743, las aulas públicas en las que se impartían primeras letras. Con la citada expulsión, la primera planta del edificio es convertida en seminario, conocido durante unos años como "Seminario Viejo", en contraposición con el nuevo Seminario de San Carlos Borromeo. Del edificio del "Seminario Viejo", es decir de las aulas donde impartía clases María Andresa, han llegado hasta nosotros sus planos de 1778

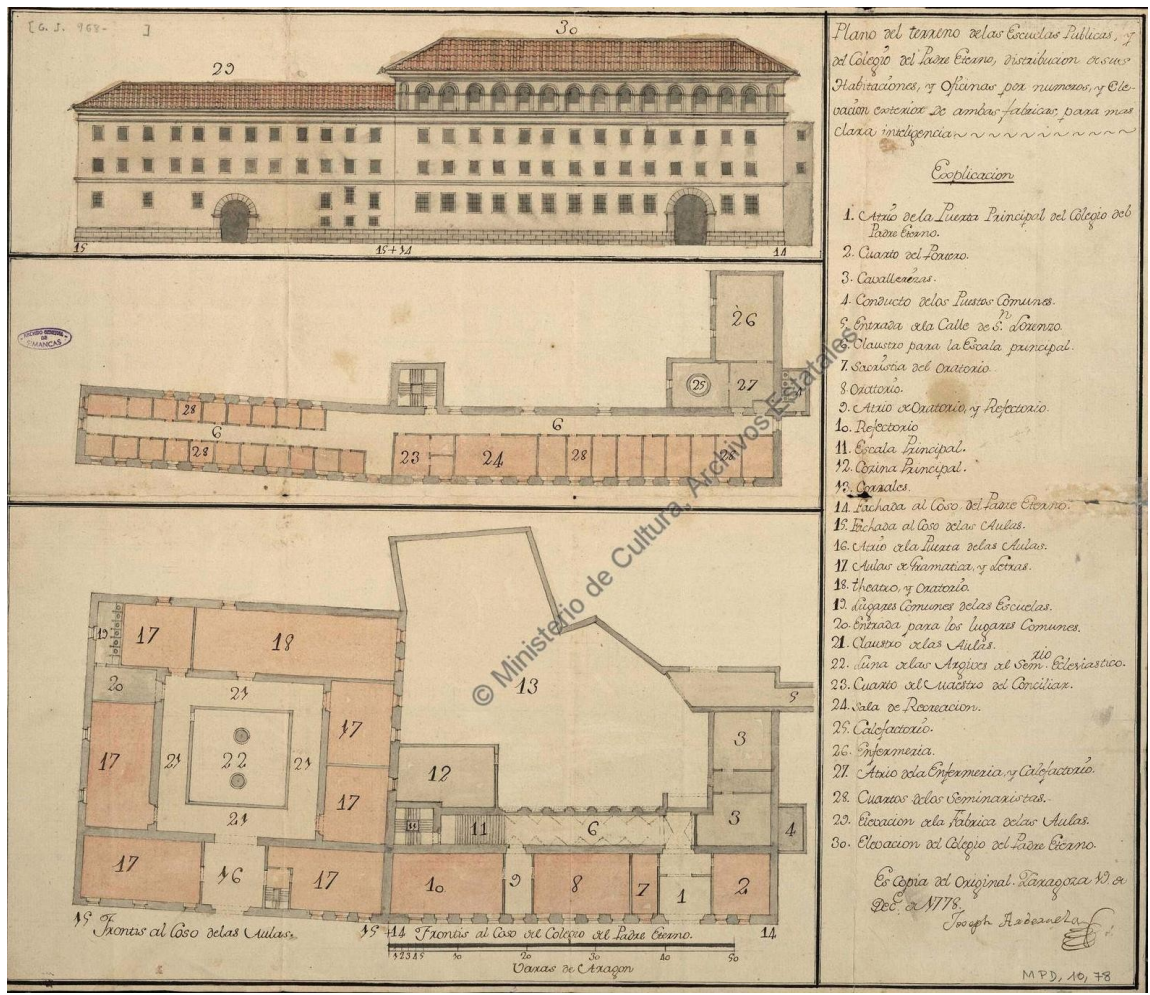

Plano de las aulas públicas donde trabajó María Andresa y un grabado que nos recuerda su uso como polvorín en los sitios de 1808, saltando por los aires en un accidente fortuito.

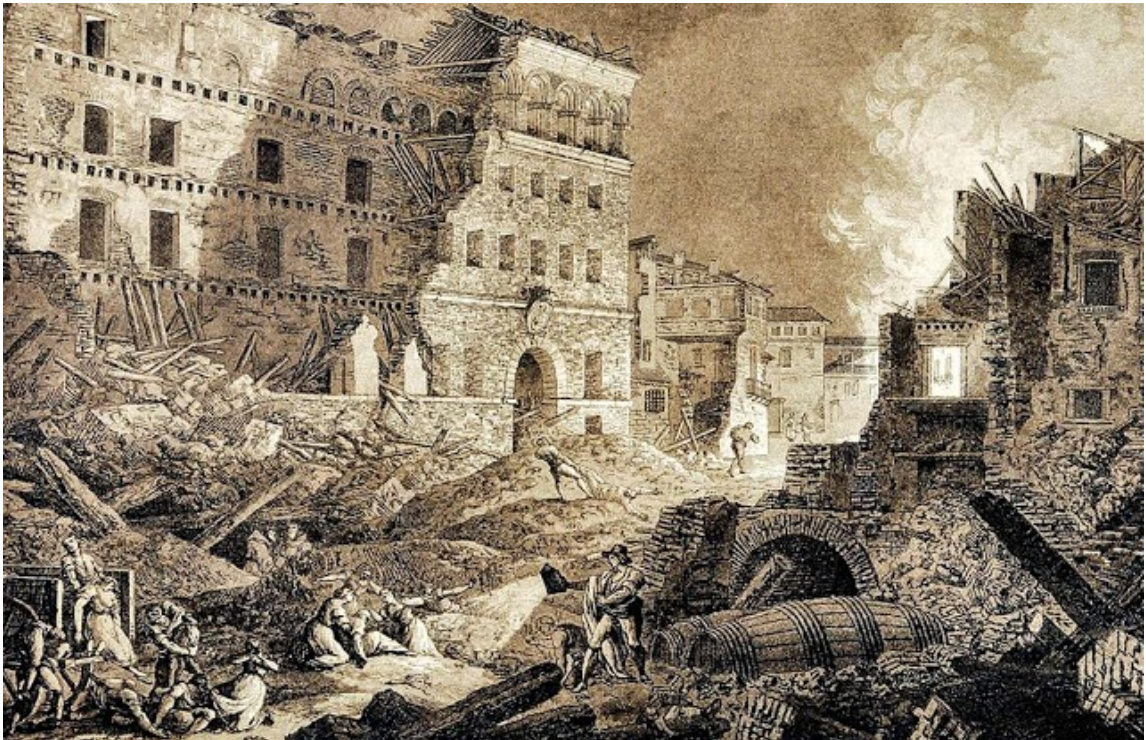
Grabado del edificio del seminario viejo después del accidente

El 23 de octubre de 1780 fallece María Casamayor. Después de recibir los sacramentos de Penitencia, Viático y Extremaunción, su cuerpo fue enterrado en el cementerio de la iglesia del Pilar [7].

**Bibliografía**

**Notas**

[1] Juan Joseph procede de Oloron (Francia) y es hijo de Juan Casamayor y de María Abales, aunque llegará solo a Zaragoza. (El apellido materno "mancebo" que aparece en algunas publicaciones recientes es por tanto erróneo). Las capitulaciones matrimoniales estipulan que él aportará al matrimonio 2.000 libras jaquesas, una respetable cantidad y la pareja deberá vivir durante más de un año en casa de los padres de élla. Juana Rosa de La Coma es hija de Juan de La Coma y de María Alexandre, vecinos de

Zaragoza. Juana Rosa tiene dos hermanos menores, Joseph, el que será heredero y Thomas, que será sacerdote del lugar de Monegrillo. *Archivo de protocolos de Zaragoza, nº* 151. Notario Blas de Villanueva**.** *p.620ss: Capitulación Matrimonial de Juan Joseph Casamayor y Juana Rosa La Coma 1-3-1705.*

[2] Libros Sacramentales del *Archivo Capitular del Pilar, Libros 4 1689-1717 y Libro 5 1718-1735.* De los 9 hijos, 2 habrán fallecido en el momento de realizar el censo de población de 1733. María Juana Rosa Andresa es bautizada el 1-12-1720 actuando como padrino su hermano Juan Casamayor y como madrina de honor su hermana Valera Martina Casamayor. Tres años más tarde el 18-4-1723 María Andresa recibirá la confirmación.

[3] La familia vivirá hasta la muerte del padre en una casa alquilada de la calle del Pilar. Al principio contarán con la ayuda del abuelo materno de María Andresa, Juan de La Coma. *Archivo de protocolos de Zaragoza, nº5557* Notario Roque Antonio Nuñez*: Arrendamiento 10-2-1708.*

[4] Joseph de La Coma, pedirá en 1740 un préstamo (un censal) de 400 libras jaquesas al Capítulo de la Iglesia de Santa Cruz, incluyendo como aval las dos casas que poseía la familia La Coma, una en la calle Albardería (parroquia de S. Pablo) y la de su residencia en la calle al horno de la Yedra (parroquia del Pilar). Al no tratarse de una gran cantidad de dinero, nos hace pensar en algún suceso familiar grave que provocará que en 1746 se inicie el proceso judicial que terminará en 1748 con el traspaso de ambas viviendas al Capítulo. *Archivo Municipal de Zaragoza. Catastro 1737. Caja Letra J. Joseph LaComa.*

[5] La familia Casamayor La Coma no tendrá propiedades (casas, campos...) y vivirán del producto de los negocios de Juan Joseph. Al morir éste desaparece la única fuente de sustento familiar. Así, veremos que al hermano de María Andresa, Juan Gregorio como eclesiástico o a su hermana Juana casada con un mercader francés. Juana Rosa, la madre, fallecerá más tarde en 1764. Muerte del padre en *Archivo Capitular del Pilar, Libro 6 1736-1756 p.389v.* Muerte de la madre en *Archivo del Pilar, Libro 7 1756-1772 p.357r*

[6] Estaba situado en la actual calle Coso entre las calles S. Jorge y S. Lorenzo.

[7] La nota de defunción (*Archivo Capitular del Pilar, Libro 8. 1773-1791),* conocida desde Latassa, termina con una dirección, *calle de La Coma*. Esta nota no hace referencia a donde vivía, sino a la dirección a la que había que dirigirse para cobrar los gastos del sepelio, quizá siguiendo indicaciones de su hermano y Beneficiado del Pilar, Juan Gregorio Casamayor, que residía en la bocacalle de la calle La Coma llamada del *Horno de la Caraza.* La calle La Coma (que iba desde la plaza del Pilar a la calle Santiago, conocida en la actualidad y desde 1867 como la calle Damián Forment) se denominaba en vida de María Andresa, calle de Subida al Horno de la Yedra. Desconocemos el motivo del cambio de nombre, aunque hay que pensar que se debió a la actividad de los comerciantes Juan de La Coma (abuelo materno de María Andresa, fallecido en 1718) o más seguramente de Joseph de La Coma (tío de Andresa y heredero de la familia) que vivieron en esa misma calle hasta la ruina de la familia en 1748.

## Agradecimientos

**Quien tiene este libro, tiene un tesoro.**

El único ejemplar conocido que se conserva del *Tyrocinio Arithmetico* se encuentra físicamente en la Biblioteca Nacional con signatura 3/49730. A ésta llega procedente de la Biblioteca Real, la cual fue fundada por Felipe V en 1711. En la portadilla y en la primera página del libro se aprecia un doble sello de la Biblioteca Real, con la iniciales **B.R**.

En la contratapa de la encuadernación hay una anotación manuscrita "Amaya in codicem" que podría referirse a una cita de una obra jurídica de Francisco de Amaya (1585-1640) sin conexión con el propio libro.

No se ha localizado ningún ejemplar más ni en bibliotecas públicas ni privadas. Manuel Jiménez Catalán en su obra *Ensayo de una tipografía Zaragoza del siglo XVIII* (1929) dedica la entrada 1516 al Tyrocinio Arithmetico, escribiendo lo siguiente.

*Tirocinio Arithmetico. Instrucción de las quatro reglas llanas, que saca a la luz Casandro Mamés de la Marca y Arioa y lo dedica a la Escuela Pía en su Colegio de Santo Tomás de Zaragoza.- Zaragoza, 1738.*

*Un volumen en 4º.*

*El autor de esta obra es Dª Maria Andrea Casamayor y de la Conca.*

*Cat. Gasca*

La librería de Cecilio Gasca (calle Coso, 33, Zaragoza) publicó a finales del siglo XIX y en torno a 1900 varios catálogos con una periodicidad irregular y desconocida. En estos folletos se daba publicidad a los libros existentes en su establecimiento. Manuel Jiménez consultaba estos catálogos para completar su obra, y gracias a ese comentario conocemos que a finales del siglo XIX al menos un ejemplar se encontraba disponible en las estanterías del librero Cecilio Gasca.

Julio Bernues, Pedro J. Miana

Departamento de Matemáticas

IUMA & Universidad de Zaragoza

bernues@unizar.es, pjmiana@unizar.es